\documentclass[12pt]{article}
\newcommand{\copyleft}{
GNU FDL\thanks{
Copyright (C) 2003, 2004, 2006 Peter G. Doyle.
Permission is granted to copy, distribute and/or modify this document
under the terms of the GNU Free Documentation License, 
as published by the Free Software Foundation;
with no Invariant Sections, no Front-Cover Texts, and no Back-Cover Texts.
}}
\title{$\Tetra$ and $\Didi$, the cosmic spectral twins}
\author{
Peter G. Doyle
\and
Juan Pablo Rossetti
\thanks{Supported by Conicet and a Guggenheim fellowship.}
}
\date{Version dated 31 May 2006\thanks{Differs only very slightly
from the version published in \emph{Geometry and Topology}
vol. 8 (2004), pages 1227--1242.}
\\ \copyleft}

\usepackage{graphics, epsfig}
\newcommand{\putfig}[3]{
\begin{figure}[cbt]
\centerline{\mbox{\includegraphics*[width=5.0truein]{figures/#2.eps}}}
\caption{#3}
\label{fig:#1}
\end{figure}
}
\newcommand{\figref}[1]{\ref{fig:#1}}

\newcommand{\defineterm}{\emph}
\newcommand{\Lap}{\Delta}
\newcommand{\dddxdx}{\frac{\partial^2}{\partial x^2}}
\newcommand{\dddydy}{\frac{\partial^2}{\partial y^2}}
\newcommand{\dddzdz}{\frac{\partial^2}{\partial z^2}}
\newcommand{\TwoTall}{\mathrm{TwoTall}}
\newcommand{\Tetra}{\mathrm{Tetra}}
\newcommand{\Didi}{\mathrm{Didi}}
\newcommand{\Z}{\mathbf{Z}}
\newcommand{\R}{\mathbf{R}}
\newcommand{\Zmodtwo}{\Z/2\Z}
\newcommand{\Zmodfour}{\Z/4\Z}
\newcommand{\mod}{/}
\newcommand{\cross}{\times}

\newcommand{\after}{\circ}
\newcommand{\half}{\frac{1}{2}}
\newcommand{\quarter}{\frac{1}{4}}
\newcommand{\linearspan}[1]{\langle #1 \rangle}
\newcommand{\sigmatetra}{\sigma^\Tetra}
\newcommand{\sigmadidi}{\sigma^\Didi}
\newcommand{\symtetra}[1]{\sigmatetra(#1)}
\newcommand{\symdidi}[1]{\sigmadidi(#1)}
\newcommand{\fig}[1]{}
\newcommand{\dee}{\mathrm{d}}

\newcommand{\Vol}{\mathrm{Vol}}
\newcommand{\choices}[1]{
\left\{
\begin{array}{ll}
#1
\end{array}
\right.
}
\begin{document}
\maketitle
\begin{abstract}
We introduce a pair of isospectral but non-isometric compact
flat 3-manifolds called $\Tetra$ (a \defineterm{tetracosm})
and $\Didi$ (a \defineterm{didicosm}).
The closed geodesics of $\Tetra$ and $\Didi$ are very different.
Where $\Tetra$ has two quarter-twisting geodesics of the
shortest length, $\Didi$ has four half-twisting geodesics.
Nevertheless, these spaces are isospectral.
This isospectrality can be proven directly by matching eigenfunctions
having the same eigenvalue.
However, the real interest of this pair---and what led us to discover
it---is the way isospectrality emerges from the Selberg trace formula,
as the result of a delicate interplay between the lengths and
twists of closed geodesics.
\end{abstract}

\subsection*{Introducing $\Tetra$ and $\Didi$}

A \defineterm{platycosm} is a compact flat 3-manifold.
Simplest among platycosms are the \defineterm{torocosms}
(the artifacts formerly known as `3-dimensional tori').
Torocosms come in various shapes and sizes.
Among these, we distinguish the
\defineterm{cubical torocosm}
$\R^3 \mod \Z^3$,
and the \defineterm{two-story torocosm}
$\TwoTall = \R^3 \mod (\Z \cross \Z \cross 2\Z)$.

All other platycosms arise as quotients of torocosms.
There are 10 distinct types in all, of which 6
(torocosm; dicosm; tricosm; tetracosm; hexacosm; didicosm)
are orientable.
The spaces themselves are well known,
but the naming scheme, due to Conway, is new.
The naming scheme and the spaces themselves are
described in great detail
by Conway and Rossetti
in \cite{conwayrossetti:seeing}.
The spaces are described under different names by Weeks
\cite{weeks:sos},
and Weeks
(see \cite{weeks:curved})
has also produced software which allows you to
`fly around' inside these spaces, and many others as well.

Here we are concerned with two specific platycosms:
$\Tetra$, a \defineterm{tetracosm}, and
$\Didi$, a \defineterm{didicosm}.
Please note that the prefix `didi-' is
a doubling of the prefix `di-', and not some exotic Greek root.
The word `didicosm' is pronounced `die-die-cosm', but
$\Didi$ is pronounced `Dee-dee'.
$\Tetra$ and $\Didi$ turn out to be, up to scale, the unique pair of
\defineterm{cosmic spectral twins} (non-isometric
platycosms with identical Laplace spectrum).

$\Tetra$ and $\Didi$ are both 4-fold quotients of $\TwoTall$.
$\Tetra$ is the quotient of $\TwoTall$ by a
fixed-point-free action of $\Zmodfour$,
while $\Didi$ is the quotient by a
fixed-point-free action of $\Zmodtwo \cross \Zmodtwo$.
To get $\Tetra$, we adjoin to the translation group
$\Z \cross \Z \cross 2 \Z$ in $(x,y,z)$-space the
quarter-turn screw motion
\[
\tau:(x,y,z) \mapsto (-y,x,z+1/2)
.
\]
To get $\Didi$, we adjoin instead the two half-turn screw motions
\[
\rho_x:(x,y,z) \mapsto (x+1/2,-y,-z)
\]
and
\[
\rho_y:(x,y,z) \mapsto (-x,y+1/2,1-z)
,
\]
which together with the translations generate a third half-turn screw motion
\[
\rho_z:(x,y,z) \mapsto (1/2-x,1/2-y,z+1)
.
\]

Both $\Tetra$ and $\Didi$ have as a fundamental domain the box
\[
[-1/2,1/2] \cross [-1/2,1/2] \cross [0,1/2]
,
\]
and in both cases the four vertical sides are glued up in parallel
in the usual way,
front to back and left to right, yielding
a stack of square tori.
The difference comes in the glueings of the top and bottom.
(See Figure \figref{10}.)
\putfig{10}{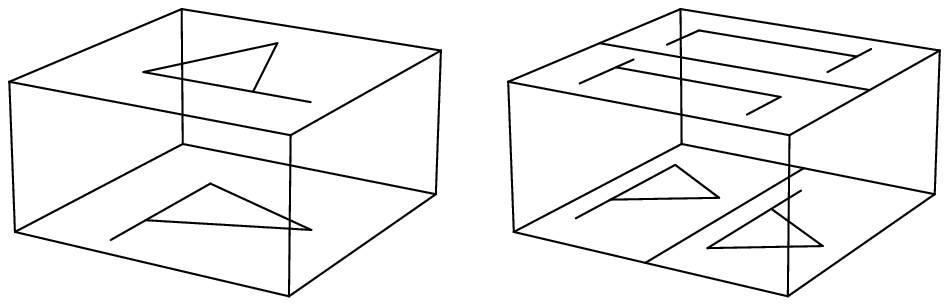}{$\Tetra$ and $\Didi$.
The sides of the box glue
back to front and left to right in the usual way;
the tops and bottoms glue as indicated.
Note that in the case of $\Didi$, the top and bottom
glue not to each other but \emph{each to itself},
yielding two Klein bottles embedded in the quotient
(which is nonetheless orientable!).}
To get $\Tetra$, you use $\tau$ to glue the bottom to the top with a
quarter-turn.
To get $\Didi$, you use $\rho_x$ to glue the bottom to itself
via a glide reflection,
and $\rho_y$ to glue the top to itself via a glide reflection.
These glueings produce two Klein bottles embedded in $\Didi$.
There is also a third, `vertical'  Klein bottle, associated to $\rho_z$.

{\bf Note.}  We have described $\Tetra$ and $\Didi$ as quotients of
a common 4-fold cover.
In fact they have a common 2-fold cover, the `cubical dicosm',
and a common 2-fold orbifold quotient.
While both of these related spaces have a role to play
in explaining the relationship between
$\Tetra$ and $\Didi$,
we won't have any further occasion to discuss them here.

\subsection*{Non-isometric}

$\Tetra$ and $\Didi$ are not isometric.
In fact, since they have different fundamental groups, they are
not even homeomorphic.
($\Tetra$ has first Betti number $1$,
while $\Didi$ has first Betti number $0$.)

Moreover, in contrast to many of the known examples of spectral
twins,
their closed geodesics are markedly different.
In a platycosm, when you go around a closed geodesic
you come back twisted through some angle $\theta$.
In a torocosm, this twist angle is always $0$.
But in $\Tetra$, the shortest geodesics
have twist $\theta = \pi/2$:
We call such geodesics \defineterm{quarter-twisting geodesics},
or \defineterm{quarter-twisters}.
In $\Didi$,
the shortest geodesics have twist $\theta=\pi$:
We call such geodesics \defineterm{half-twisting geodesics},
or \defineterm{half-twisters}.
(See Figure \figref{20}.)

The fact that the shortest geodesics in $\Tetra$ are
quarter-twisters while those in $\Didi$ are half-twisters
already shows that these spaces are non-isometric.
(Indeed, $\Didi$ has no quarter-twisting geodesics at all;
this is an aspect of the fact that $\Tetra$ and $\Didi$ have different
`holonomy groups', namely $\Zmodfour$ versus $\Zmodtwo \cross \Zmodtwo$.)

Let us count the short geodesics in $\Tetra$ and $\Didi$.

{\bf Warning.}
When counting geodesics,
we count each pair of oppositely-oriented geodesics only once.

\putfig{20}{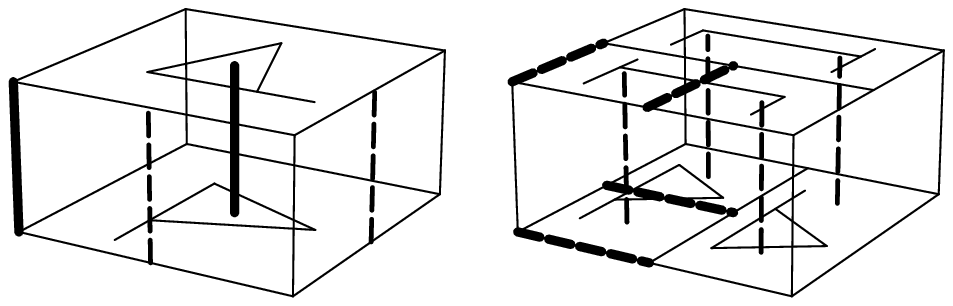}{Twisted geodesics in $\Tetra$ and $\Didi$.
Solid thick segments close up into quarter-twisting geodesics
of length $1/2$;
dashed thick segments close up into half-twisting geodesics
of length $1/2$;
dashed thin segments join together in pairs to form
half-twisting geodesics of length $1$.
Note that in the case of $\Didi$, if you start up one of the dashed thin
vertical segments, you continue down the segment `kitty-corner' to it.}

In $\Tetra$ there are two quarter-twisting geodesics of length $1/2$,
one running up the middle of the box along the line $x=y=0$,
and one running up the four identified edges of the box.
The vertical midlines of the four sides of the box combine to
give a third geodesic, but this one is a half-twisting geodesic
of length 1.

In $\Didi$, there are four half-twisting geodesics of length $1/2$,
two associated with $\rho_x$ sitting
in the Klein bottle gotten by glueing the bottom of the box, and two
associated with $\rho_y$ sitting in the Klein bottle gotten by glueing
the top.
(See Figure \figref{20}.)
In addition, there are two half-twisting geodesics sitting in the
$\rho_z$ Klein bottle, but these have length $1$.

Let's call a geodesic \defineterm{twisted} if it has a nontrivial twist.
We've located $3$ twisted geodesics in $\Tetra$,
and $6$ in $\Didi$.
These are all \defineterm{primitive}, which means that they don't
arise by going more than once around a shorter geodesic.
In fact these are the only primitive twisted geodesics in $\Tetra$ and 
$\Didi$.
Of course there are also \defineterm{imprimitive} twisted geodesics, which come from
going around a half-twister an odd number of times, or around a
quarter-twister a number of times not divisible by $4$.

{\bf Note.} To verify that we have identified all the twisted geodesics,
observe that any twisted geodesic
in $\Tetra$ (or $\Didi$) unwraps to a straight line in $R^3$.
Translating-with-a-twist along this straight line will
be a `covering transformation'---that is, one of the symmetries
of $\R^3$ consistent with all patterns obtained by unwrapping
patterns in the quotient $\Tetra$ (or $\Didi$).
The translation is by the length of the geodesic,
and the twist is equal (and opposite) to the twist of the geodesic.
In the case of $\Tetra$, the line must be vertical, because
all covering translations-with-a-nontrivial-twist run vertically.
In the case of $\Didi$, the line can run in any of the three coordinate
directions.
Look carefully at the possibilities, and you'll
see that in listing twisted geodesics
we've accounted for all of them.

\subsection*{Isospectral}

While $\Tetra$ and $\Didi$ are not isometric, they are \defineterm{isospectral}.
By definition, two spaces are isospectral if there exists some
way of matching up the eigenfunctions of the Laplacian of the two spaces
so that corresponding eigenfunctions have the same eigenvalue.
In this section,
we will show that $\Tetra$ and $\Didi$ are isospectral by describing
such a correspondence.

We are giving this explicit proof because it is entirely elementary---it
relies only on linear algebra and Fourier series---and
because it is illuminating in its own way.
Other, more `conceptual' proofs are available.
Further along, we will outline one such proof,
by way of the Selberg trace formula.
A third proof can be obtained using the
general machinery for flat manifolds developed by Miatello and Rossetti
in
\cite{mr:p},
and a fourth using the `dual' approach of
\cite{mr:flat}.
A close relative of this fourth proof
emerges naturally in the proof that $\Tetra$ and $\Didi$ are
the unique pair of cosmic spectral twins 
\cite{conwayrossetti:hearing},
discussed briefly below.

A function on $\Tetra$ corresponds to a function $f$ on $\TwoTall$
that is invariant under $\tau$, in that
\[
f = f \after \tau
.
\]
Given any function $f$ on $\TwoTall$, we can symmetrize under $\tau$ to get
a $\tau$-invariant function
\[
\symtetra{f} = \frac{1}{4}(f + f \after \tau + f \after \tau \after \tau
+ f \after \tau \after \tau \after \tau)
,
\]
which we think of as a function on $\Tetra$.
Similarly, we can get functions on $\Didi$ via the symmetrization
\[
\symdidi{f} = \frac{1}{4}(f + f \after \rho_x + f \after \rho_y + f \after \rho_z)
.
\]

Now any function $f$ on $\TwoTall$ can be written as a Fourier series:
\[
f(x,y,z) = \sum_{(a,b,c) \in \Z \cross \Z \cross \half \Z}
\hat{f}(a,b,c) \exp(2 \pi i (ax+by+cz))
.
\]
Note that the sum runs over the lattice $\Z \cross \Z \cross \half \Z$,
which is `dual' to the original lattice $\Z \cross \Z \cross 2 \Z$:
The frequencies $a$ and $b$ in the $x$ and $y$ directions are integers,
but the frequency $c$ in the $z$ direction is allowed to be a half-integer,
because the scale of the lattice in that direction is twice the scale
in the $x$ and $y$ directions.

The Fourier basis functions $\phi_{a,b,c} = \exp(2 \pi i (ax+by+cz))$,
$(a,b,c) \in \Z \cross \Z \cross \half \Z$ are eigenfunctions of
the (positive) Laplacian $\Lap = -(\dddxdx + \dddydy + \dddzdz)$:
\[
\Lap \phi_{a,b,c} = 4 \pi^2 (a^2 + b^2 + c^2) \phi_{a,b,c}
.
\]

Symmetrizing these Fourier basis functions under $\tau$ yields a spanning set
$\symtetra{\phi_{a,b,c}}$ for the functions on $\Tetra$.
This spanning set is far from being a basis.
For one thing, symmetrization lumps the basis functions
together in groups, generally of size four.
More important,
symmetrizing a basis function can kill it off altogether.
For example, 
\[
\symtetra{\phi_{0,0,1/2}}=
\symtetra{\phi_{0,0,1}}=
\symtetra{\phi_{0,0,3/2}}=
0
.
\]
However, by eliminating such redundancies,
we can prune down to a basis of
(unnormalized) eigenfunctions $\symtetra{\phi_{a_i,b_i,c_i}}$ on $\Tetra$,
with corresponding eigenvalues
$4 \pi^2(a_i^2 + b_i^2 + c_i^2)$.

\newcommand{\vary}[1]{{#1^\prime}}
Similarly, we can get a basis of eigenfunctions
$\symdidi{\phi_{\vary{a}_i,\vary{b}_i,\vary{c}_i}}$ on $\Didi$,
with corresponding eigenvalues
$4 \pi^2(\vary{a}_i^2 + \vary{b}_i^2 + \vary{c}_i^2)$.
If we can arrange that
\[
a_i^2 + b_i^2 + c_i^2 = \vary{a}_i^2 + \vary{b}_i^2 + \vary{c}_i^2
\]
for all $i$, then we will have verified that
$\Tetra$ and $\Didi$ are spectral twins.

To find such a correspondence, we will take advantage of the fact
that our two symmetrization mappings lump the basis functions
together in two different but very nearly compatible ways.
This leads to a correspondence between eigenfunctions which
is for the most part very straight-forward.
There are only two exceptional cases
that must be treated carefully.

Let
\[
V_{a,b,c}
=
\linearspan{\phi_{\pm a, \pm b, \pm c}, \phi_{\pm b, \pm a, \pm c}}
,
\]
where we're using angle brackets to denote linear span.
Please observe that $V_{a,b,c}=V_{b,a,c}$, $V_{a,b,c}=V_{-a,b,c}$, etc.
In the generic case,
namely when $a,b,c\neq 0$ and $|a| \neq |b|$,
the vector space $V_{a,b,c}$ is 16-dimensional,
and both $\sigmatetra$ and $\sigmadidi$ lump the 16 basis functions
together in groups of 4.
In this case $\sigmatetra(V_{a,b,c})$ and $\sigmadidi(V_{a,b,c})$
are both 4-dimensional,
and we can clearly take bases of these spaces
and match them up.

If it were true that
$\dim \sigmatetra(V_{a,b,c}) = \dim \sigmadidi(V_{a,b,c})$,
for all $(a,b,c) \in \Z \cross \Z \cross \half \Z$,
we would be all set.
In fact this equality holds as long as no two of the
parameters $a,b,c$ vanish,
because in these cases
$\sigmatetra$ and $\sigmadidi$ continue to lump the basis functions
together in groups of 4.
Of course $\dim \sigmatetra(V_{0,0,0}) = \dim \sigmadidi(V_{0,0,0})=1$,
so that case is no problem.
And if $c$ is a half integer, then
$\dim \sigmatetra(V_{0,0,c}) = \dim \sigmadidi(V_{0,0,c})=0$.

So the question comes down to how to handle the cases $V_{n,0,0}$ and
$V_{0,0,n}$, with $n$ a non-zero integer, which we may assume is positive.
(Remember that negating any of $a,b,c$ does not change the space $V_{a,b,c}$.)
To extend the correspondence between eigenfunctions,
we must take these remaining exceptional cases in combination.
Here is how it goes.

{\bf Odd exceptional case.}
When $n$ is a positive odd integer, 
\begin{eqnarray*}
\symtetra{V_{n,0,0}} &=&
\linearspan{\cos 2 \pi n x + \cos 2 \pi n y }
;\\
\symtetra{V_{0,0,n}} &=&
0
;\\
\symdidi{V_{n,0,0}} &=&
0
;\\
\symdidi{V_{0,0,n}} &=&
\linearspan{ \cos 2 \pi n z }
.
\end{eqnarray*}
Taken together, these cases contribute a single eigenfunction of
eigenvalue $4 \pi^2 \cdot n^2$ to the spectra of both $\Tetra$ and $\Didi$.

{\bf Even exceptional case.}
When $n$ is a positive even integer,
\begin{eqnarray*}
\symtetra{V_{n,0,0}} &=& 
\linearspan{ \cos 2 \pi n x + \cos 2 \pi n y }
;\\
\symtetra{V_{0,0,n}} &=&
\linearspan{ \exp (2 \pi i n z), \exp (- 2 \pi i n z) }
;\\
\symdidi{V_{n,0,0}} &=&
\linearspan{ \cos 2 \pi n x , \cos 2 \pi n y }
;\\
\symdidi{V_{0,0,n}} &=&
\linearspan{ \cos 2 \pi n z }
.
\end{eqnarray*}
Taken together, these cases contribute three independent eigenfunctions
of eigenvalue $4 \pi^2 \cdot n^2$ to both spectra.

By matching up these exceptional cases as indicated,
we finish the job of matching up eigenfunctions of $\Tetra$ and $\Didi$,
and thus concretely demonstrate that these spaces
are spectral twins.

Note that our scheme for matching eigenfunctions
involves some arbitrary, symmetry-breaking choices.
This shows up clearly in the even exceptional case above,
but it is an issue even in the `generic' case.
We will see this same kind of symmetry-breaking
again when we look at the proof of
isospectrality by way of the Selberg trace formula.

\subsection*{Unique}

$\Tetra$ and $\Didi$ are, up to scale, the \emph{only}
pair of non-isometric isospectral platycosms:
They are the \emph{two-and-only} cosmic spectral twins.
The proof, due to Rossetti,
involves a case-by-case analysis of all possible
spectral coincidences among and between platycosms of
the 10 possible types.
In \cite{conwayrossetti:hearing},
Rossetti and Conway give a streamlined version of this proof,
using Conway's theory of lattice conorms as an
organizing principle.
As you would expect,
the techniques used in proving uniqueness
yield another proof that $\Tetra$ and $\Didi$ are spectral twins.

A key ingredient in the uniqueness proof is Schiemann's theorem
\cite{schiemann:three}
that there are no spectral twins among torocosms:
If $\R^3/\Lambda_1$ and $\R^3/\Lambda_2$ are isospectral, then
they (and the lattices $\Lambda_1$ and $\Lambda_2$) are isometric.
Milnor's original example
of spectral twins was a pair of
16-dimensional tori
\cite{milnor:sixteen}.
Subsequently, lower-dimensional pairs of isospectral tori were found,
culminating with the discovery of a 4-dimensional pair by Schiemann,
simplified and extended to a 4-parameter family of pairs by Conway and Sloane
\cite{conwaysloane:four}.
Schiemann showed that as far as tori are concerned,
dimension 4 is the end of the line.
By opening the field up to other flat manifolds,
we can get down to dimension 3---but just barely!

\subsection*{Selberg}

Here, as promised above, we outline a proof of isospectrality
by way of the Selberg trace formula.
The version of the trace formula that we want to use expresses the
Laplace transform of the spectrum (or properly speaking,
of the spectral measure) as the sum of 
contributions attributable to families of closed geodesics.  
To show that $\Tetra$ and $\Didi$ are isospectral, we will examine the
closed geodesics of each, and check that the total spectral contribution of
$\Tetra$'s geodesics is just the same as that of $\Didi$'s.

The relevant computations are indicated in
Table \ref{table:balance}.
Here we will explain informally what lies behind the
computations in the table.
The discussion is contrived
in such a way as to allow us to put off
actually writing down Selberg's formula until after we have put it to use.
Our reason for preferring this inverted approach is that
(in the present case, at least)
the Selberg formula is easier to apply than to state.

\newcommand{\wt}[3]{#1 \cdot \frac{#2}{#3}}
\newcommand{\nl}{\\[0.1cm]}
\begin{table}
\[
\begin{array}{cc||cccc|cccc||cccc|cccc}
\multicolumn{2}{c}{\ }&
\multicolumn{8}{c}{\Tetra}&
\multicolumn{8}{c}{\Didi}
\nl
l&w_l&
n&t&k&w&
n&t&k&w&
n&t&k&w&
n&t&k&w
\nl
\hline
\half&4&
2&\quarter&1&\wt{2}{2}{1}&
\!\!\phantom{8^{8^8}}\!\!
&&&&
4&\half&1&\wt{4}{1}{1}
&&&&
\nl
1&2&
2&\half&2&\wt{2}{1}{2}&
1&\half&1&\wt{1}{1}{1}&
&&&&
2&\half&1&\wt{2}{1}{1}
\nl
\frac{3}{2}&\frac{4}{3}&
2&\quarter&3&\wt{2}{2}{3}&
&&&&
4&\half&3&\wt{4}{1}{3}
&&&&
\nl
2&0
&&&&
&&&&
&&&&
&&&&
\nl
\frac{5}{2}&\frac{4}{5}&
2&\quarter&5&\wt{2}{2}{5}&
&&&&
4&\half&5&\wt{4}{1}{5}
&&&&
\nl
3&\frac{2}{3}&
2&\half&6&\wt{2}{1}{6}&
1&\half&3&\wt{1}{1}{3}&
&&&&
2&\half&3&\wt{2}{1}{3}
\nl
\frac{7}{2}&\frac{4}{7}&
2&\quarter&7&\wt{2}{2}{7}&
&&&&
4&\half&7&\wt{4}{1}{7}
&&&&
\nl
4&0
&&&&
&&&&
&&&&
&&&&
\nl
\frac{9}{2}&\frac{4}{9}&
2&\quarter&9&\wt{2}{2}{9}&
&&&&
4&\half&9&\wt{4}{1}{9}
&&&&
\\
\multicolumn{2}{c||}{\ldots}&
\multicolumn{4}{c|}{\ldots}&
\multicolumn{4}{c||}{\ldots}&
\multicolumn{4}{c|}{\ldots}&
\multicolumn{4}{c}{\ldots}
\end{array}
\]
\caption{Balancing geodesics.
This table shows the balancing of the spectral contributions from
the twisted geodesics in $\Tetra$ and $\Didi$.
Here $l$ is length, 
and $w_l$ the total spectral contribution (\defineterm{weight})
of geodesics of length $l$,
measured in units of the spectral contribution
of a primitive half-twisting geodesic of length $l$.
The point of this table is to demonstrate that
$w_l$ is the same for $\Tetra$ and $\Didi$.
For geodesics of a specific kind,
$n$ tells the number of geodesics;
$t$ the twist (either $\quarter$ or $\half$);
$k$ the imprimitivity exponent;
and $w$ the aggregate spectral weight
for geodesics of this kind.
An individual geodesic with imprimitivity exponent $k$
has weight $1/k$ if it is half-twisting,
and $2/k$ if it is quarter-twisting.
Weights do not depend on the handedness of the twist,
so we do not distinguish
between $1/4$-twisting and $3/4$-twisting geodesics.
}
\label{table:balance}
\end{table}

Recall that when it comes to the shortest
geodesics, which have length $1/2$, $\Tetra$ has two `quarter-twisters',
while $\Didi$ has four `half-twisters'.
Now it happens that, in a flat 3-manifold, the spectral contribution
of any primitive quarter-twisting geodesic is just twice that of
a primitive half-twisting geodesic.
(This is an aspect of a general phenomenon:
`The more the twist; the less the contribution.')
So as far as the shortest geodesics
go, the contributions to the spectrum are the same.

Next come geodesics of length 1.
Both $\Tetra$ and $\Didi$ have two 2-dimensional families of non-twisting
geodesics, which they inherit from the common cover $\TwoTall$.
These common families of non-twisting geodesics
make identical contributions to the spectrum.
In general, the non-twisting geodesics in $\Tetra$ and $\Didi$ are all inherited
from the common cover,
and consequently contribute equally to their spectra.
So we don't have to worry about non-twisting geodesics.

Looking at twisted geodesics of length 1,
the only kind that arise are half-twisters.
In $\Tetra$ we already identified
one primitive half-twister running vertically up the midlines of the
sides of the box we have chosen as our fundamental domain;
in $\Didi$, we have two primitive half-twisters sitting in the vertical Klein
bottle.
That's it, as far as primitive geodesics are concerned.
However, in $\Tetra$, we also have two imprimitive half-twisters,
gotten by running twice around
those two primitive quarter-twisters of length $1/2$.
In the Selberg formula the spectral contribution of an
imprimitive geodesic
must be divided by its \defineterm{degree of imprimitivity}
or \defineterm{exponent},
which is the number of times it runs around its primitive ancestor.
So the spectral contribution of $\Tetra$'s
one new half-twister and two recycled quarter-twisters
just matches that of $\Didi$'s 
two brand new half-twisters.

Next among twisted geodesics are those of length $3/2$.
Here we are back to balancing $\Tetra$'s two quarter-twisters, now
thrice-imprimitive, against $\Didi$'s four half-twisters, also
thrice-imprimitive.

At length $2$, there are no non-twisting geodesics.

Length $5/2$ is like $1/2$ and $3/2$:
$\Tetra$ has two quarter-twisters
and $\Didi$ four half-twisters, all now five-times-imprimitive

Length $3$ is like length $1$:
We are back to balancing
$\Didi$'s two half-twisters, now thrice-imprimitive,
against $\Tetra$'s one half-twister, now thrice-imprimitive,
and two quarter-twisters, now recycled as six-times-imprimitive
half-twisters.

And so it goes on up the line.  Thus $\Tetra$ and $\Didi$ are isospectral.

\subsection*{Formula}

We have chosen to describe the geodesic balancing act between
$\Tetra$ and $\Didi$ in words, rather than symbols, in order
to put off having to state explicitly Selberg's formula for platycosms.
The goal was to show how you can use the formula without having to
know precisely what the formula is.  
Now, here comes the formula.

Let $M=\R^3/\Gamma$ be a platycosm with covering group $\Gamma$,
and let $\Lambda \subset \Gamma$ be the lattice subgroup of $\Gamma$.
Let $0=\lambda_0<\lambda_1 \leq \lambda_2 \leq \lambda_3 \leq \ldots$
be the sequence of eigenvalues of the Laplacian on $M$,
where as usual multiple eigenvalues are listed multiple times
in the sequence.
For a geodesic $g$, let $l(g)$ be its length, $\theta(g)$ its
twist (in radians), and $k(g)$ its imprimitivity exponent.
Let $G$ be the set of (nontrivially) twisted closed
geodesics of $M$.
Here again, we take no notice of orientation:  oppositely-oriented geodesics
are considered to be the same.

The version of the Selberg trace formula we need
relates the spectrum of $M$ to the geometry of $M$,
by giving two separate ways to compute a certain
function $K(t)$, the `trace of the heat kernel'.
The first way is in terms of the spectrum, while
the second way is in terms of the geometry.
For present purposes, we don't need to know
just what $K(t)$ is---just that the two expressions are equal.

Here's the formula:
\[
K(t) = 
\sum_n e^{-\lambda_n t} =
\int_0^\infty \frac{1}{(4 \pi t)^{\frac{3}{2}}} e^{-\frac{s^2}{4t}} \dee N(s)
,
\]
where
\[
N(s) = \Vol(M) \left | \{\gamma \in \Lambda: |\gamma| \leq s \} \right | 
+
2 \sum_{g \in G} \frac{1}{k(g)} V_{l(g),\theta(g)}(s)
,
\]
where
\[
V_{h,\theta}(s) =
\choices{
0, & 0 \leq s < h \\
h \pi \frac{s^2-h^2}{(2 \sin \frac{\theta}{2})^2}, & h \leq s 
}
.
\]
Here $V_{h,\theta}(s)$ is the volume of a cylinder of height $h$ and
\defineterm{$\theta$-twisted height} $s$,
which we define as follows:
If the sides of a cylinder of height $h$ are replaced with parallel segments,
and the top is twisted through an angle $\theta$ relative to the bottom,
the segments stretch
to form (part of) a hyperboloid of one sheet; their stretched length is
what we're calling the $\theta$-twisted height.
(See Figure \figref{30}.)

Note the factor of $2$ in front of the sum over $G$ in the formula
for $N(s)$.
This arises because of our lumping together oppositely-oriented geodesics.
When it comes to doing explicit
computations using this formula,
this lumping appears unnatural,
because translating one way around a geodesic is
not at all the same thing as translating the other way.
But for us to try to enforce a distinction between oppositely-oriented
geodesics
would only cause confusion,
because it would run counter to established practice.
Nobody is going to be comfortable with the notion
that $\Tetra$ has \emph{four} geodesics of length $\half$.

Now, what makes Selberg's formula so useful
is the fact that the function
$K(t) = \sum_n e^{-\lambda_n t}$ determines the spectrum:
In fact, it is the Laplace transform of a mass distribution with a unit
mass placed at every eigenvalue, so we can recover the spectrum
by inverting a Laplace transform.
The trace formula thus shows that the Laplace spectrum is determined
by the function $N(s)$, which is computed directly from geometrical data.

The use we made above of the trace formula depended essentially
on only two of its properties.
First is the amazing fact that the contribution
of a twisted geodesic to $N(s)$ depends on the twist $\theta$
only through the factor $\frac{1}{\sin^2 \frac{\theta}{2}}$.
For a half-twisting geodesic ($\theta = \pi$), this factor is $1$;
for a quarter-twisting geodesic ($\theta = \pi/2$), this factor is $2$.
Second is the way recycled geodesics get only partial credit in $N(s)$,
because of the
factor $\frac{1}{k(g)}$.

Please observe that,
since we are only using this formula
to show that the two spaces $\Tetra$ and $\Didi$ are isospectral,
it isn't crucial that we have the formula exactly right.
If the formula as written
were off by a constant factor here or there,
say a stray factor of $2$ or $\sqrt{\pi}$,
it would do us no harm,
as long as the error applied equally to $\Tetra$ and $\Didi$.
That's good, of course.
Working only with `higher-level' properties of the formula
means that we're approaching things in a `conceptual' way,
and insulating ourselves from possible bugs in the formula.
The drawback is that
we haven't given the formula a real workout.

{\bf Exercise.}
Give the formula a real workout by using it to compute explicitly
the heat trace for $\TwoTall$, and then $\Tetra$.
(No need to check $\Didi$,
since we've already checked that the answer will be the same 
for $\Didi$ as for $\Tetra$.)
Observe that
\[
K_{\Tetra} - \frac{1}{4} K_{\TwoTall}
=
K_{\R/\half\Z} - \frac{1}{4} K_{\R/2\Z}
.
\]
Explain this remarkable `effective 1-dimensionality' of the relationship
between $\Tetra$ and $\TwoTall$.

\subsection*{Derivation}
So that's the Selberg trace formula for platycosms,
or rather,
one particular form of it.
As for the derivation,
equivalent formulas are derived in the papers 
of Gangolli \cite{gangolli:selberg} and Berard-Bergery \cite{bb:selberg}.
The formula can also be derived by adapting
the classical Poisson summation formula,
which was the original inspiration for the Selberg formula:
This is the approach of Sunada \cite{sunada:flat} and
Miatello and Rossetti \cite{mr:flat}.

More mundanely,
the formula can also be written down by identifying the
Laplace transform of the spectrum with the trace of the
heat kernel in the usual Selbergian way, and
evaluating the diagonal of the heat kernel
using the `method of images' from sophomore physics.
Pairs of images lying within a given distance of each other
fall into cylinders of given `twisted height'.
Measuring the volume of these cylinders requires only
the Pythagorean theorem and the usual formula for the volume of a cylinder.
The fact that a quarter-twisting geodesic contributes twice as much
to the spectrum as a half-twisting geodesic boils down to the
fact that a cylinder of given height $h$ and
quarter-twisted height $s$ has twice the volume of a cylinder
with the same $h$ and half-twisted height $s$.
(See Figure \figref{30}.)
\putfig{30}{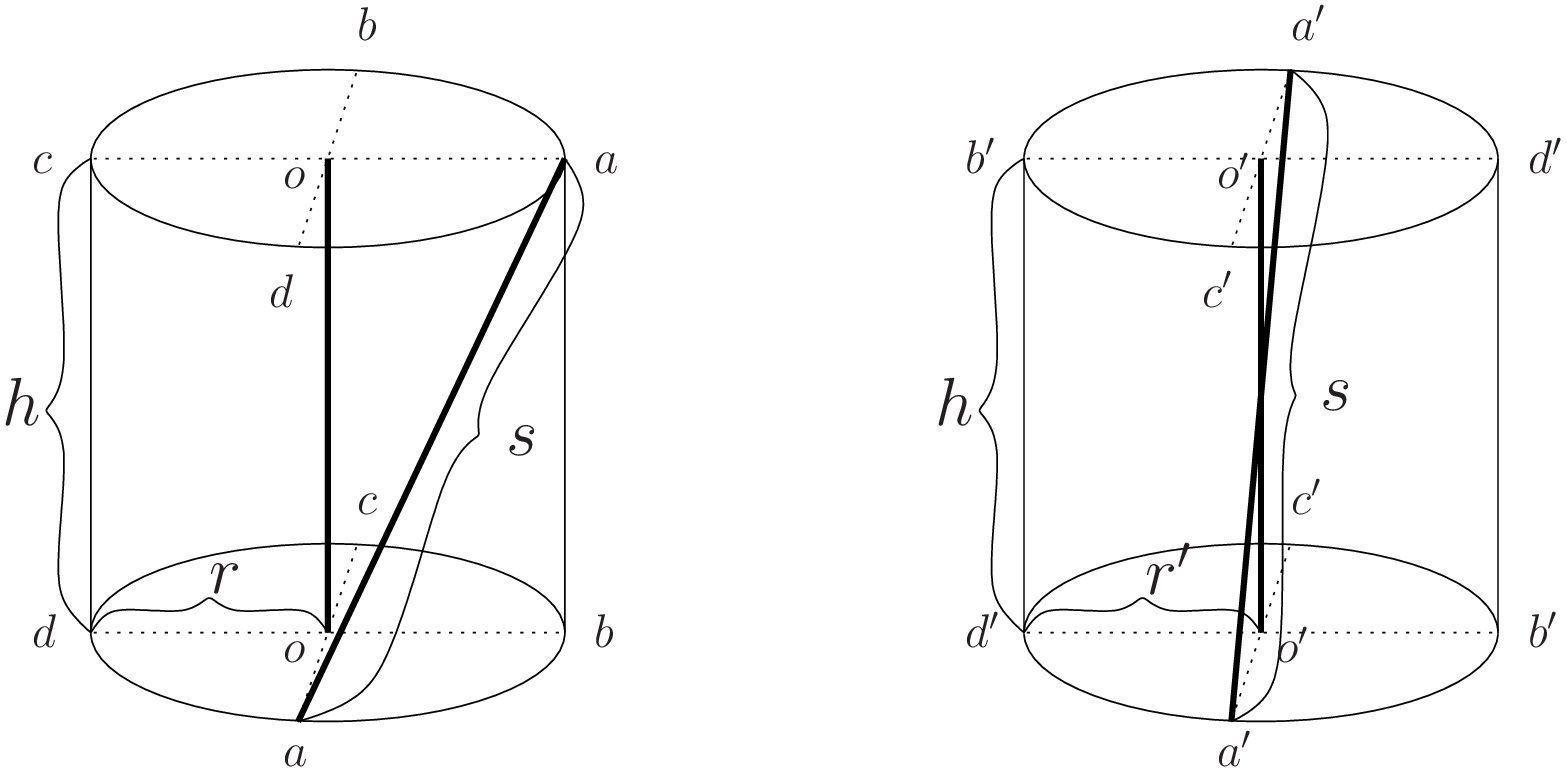}{
Volume of cylinders with given
height $h$ and twisted height $s$.
$\Vol_\quarter = \pi r^2 h$;
$\Vol_\half = \pi r'^2 h$.
By the Pythagorean theorem,
$h^2+2r^2=s^2$ and $h^2+4r'^2=s^2$.
So $\Vol_\quarter = 2 \Vol_\half$.
}

A similar hands-on approach to Selberg's formula
works for hyperbolic 2- and 3-manifolds (or orbifolds),
with the hyperbolic law of cosines taking the place of the
Pythagorean theorem.

\subsection*{Subtle isospectrality}

The remarkable Selbergian interplay between the geodesics
in $\Tetra$ and $\Didi$
is what we were looking for when we discovered this pair.
Originally, we were interested in finding (or ruling out)
an analogous pair of \emph{hyperbolic} 3-manifolds.
There are plenty of examples of spectral twins among hyperbolic
3-manifolds, but the standard methods for producing spectral twins
yield pairs whose geodesics have matching lengths and twists.
Selberg's formula seems to allow the possibility of twins that
are \defineterm{subtly isospectral},
meaning, `isospectral, but not merely by virtue of having geodesics with
matching lengths and twists'.

After some fruitless attempts to find such a pair among hyperbolic
3-manifolds,
we tried looking among flat 3-manifolds instead.

Now, the precise definition of `subtle isospectrality' for flat manifolds
is a subtle business,
because various of the geodesics in flat manifolds come in parallel families,
and you have to
be careful to choose the right definition for the `multiplicity' of
such geodesics,
which you will need in order to decide whether a give pair of manifolds
have geodesics with `matching' lengths and twists.
Of course, you could make any definition you want, and then
investigate it.
But to be interesting, the definition of multiplicity
should end up meaning that two manifolds are subtly isospectral
just if they are isospectral, but not merely because of some
straight-forward way of matching contributions to Selberg's formula.

We don't need to fuss about this here, because
however you decide to measure multiplicity,
$\Tetra$ and $\Didi$ will have different
multiplicities of geodesics of the shortest length.
That's because these shortest geodesics are all isolated,
on account of their twisting.
Thus $\Tetra$ and $\Didi$ are subtly isospectral---and
no other pairs of non-isometric platycosms are isospectral,
whether subtly or not.

The possibility of subtly isospectral hyperbolic 3-manifolds remains open.

\subsection*{Remarks}

\paragraph*{Misconception about Selberg for hyperbolic 3-manifolds.}
Some people mistakenly believe
that, for a hyperbolic 3-manifold,
Selberg's formula allows you to read off from the spectrum
the lengths and twists of the closed geodesics.
For example,
Reid
\cite{reid:isospec}
cites Gangolli \cite{gangolli:selberg}
and Berard-Bergery \cite{bb:selberg}
in support of this assertion.
Neither Gangolli nor Berard-Bergery
makes such a statement,
and Gangolli in particular is explicit
about the fact that the Selberg formula
leaves open
the possibility of an example of the kind we were (and still are)
looking for.

This misconception most likely stems from conflicting
uses of the term \defineterm{length spectrum},
which we have been studiously avoiding here.

\paragraph*{Hyperbolic surfaces.}
In another attempt to work up to hyperbolic 3-manifolds,
we looked at hyperbolic surfaces.
We were disappointed to find
\cite{doylerossetti:surfaces}
that no subtly isospectral pairs exist among hyperbolic surfaces.
The question only becomes interesting in the case of non-orientable
surfaces, because the Selberg formula immediately implies that
for hyperbolic surfaces, spectral twins always have matching lengths.
For non-orientable surfaces, it is possible to construct a plausible
scenario for matching the contributions of the geodesics of a pair of surfaces
in a way similar to that of $\Tetra$ and $\Didi$,
where geodesics don't have matching lengths.
But in the end it proves impossible to make this work:
We show that to balance the spectral contributions of geodesics on
a pair of surfaces whose
geodesics can't be matched so as to preserve length and
orientability class, you would frequently need the number of geodesics 
with length exactly $l$ to be $\geq C e^l/l$ for $C>0$.
But as a consequence of results of
Huber \cite{huber:selbergI,huber:selbergII,huber:selbergIIA}
and others, this number must be $o(e^l/l)$.
(Huber's results are stated only for orientable surfaces,
but they hold as well in the non-orientable case.)
This rules out the possibility of subtly isospectral hyperbolic
surfaces:  Even disconnected surfaces, which turn out to
be more interesting than one might imagine.

\paragraph*{Laplacian on forms.}
While $\Tetra$ and $\Didi$ are isospectral for the usual Laplacian acting on
\emph{functions},
they are not isospectral for the Laplacian acting on 1-forms or 2-forms.
This is a simple consequence of the techniques of Miatello and Rossetti
(see Theorem 3.1 of \cite{mr:p}).
However, this can also be seen immediately
because the Betti number $b_k$ equals
the multiplicity of the eigenvalue $0$ of the
Laplacian acting on $k$-forms.
$\Tetra$ and $\Didi$ have distinct $b_1$,
and hence distinct spectrum on 1-forms.
Since they are both orientable 3-manifolds, 
they have $b_1=b_2$,
so they have distinct $b_2$ and distinct
spectrum on 2-forms.
Examples of manifolds which are isospectral but have distinct
first Betti numbers, or
which are isospectral on functions but not on 1-forms,
were previously known only in higher dimensions
(see \cite{gordon:survey} and the references therein).

\bibliography{cosmic}
\bibliographystyle{hplain}

\end{document}